\documentclass[11pt]{article}
\def\ti{\tilde}

\def\0{\bf 0}
\def\1{\bf 1}
\def\t{\theta}

\def\wh{\widehat}

\def\Pb{{\bf P}}
\def\xb{{\bf x}}
\def\s{\sigma}

\def\Cov{{\bf Cov}}
\def\Eb{{\bf E}}
\def\E{{\bf E}}

\def\iy{\infty}

\def\e{\varepsilon}

\def\endsymbol{$\sqcup\mkern-12mu\sqcap$}
\def\done{\ \endsymbol\medskip}

\def\Var{\mathop{\rm\bf Var}\nolimits}
\usepackage{amsmath, amsfonts, amssymb, amsthm}

\begin{document}

\title{Adaptive variable selection in nonparametric sparse additive models}
\author{Cristina Butucea$^{1,2}$ and Natalia Stepanova$^3$\\
{\small $^1$ Universit\'e Paris-Est Marne-la-Vall\'ee, }\\
{\small LAMA(UMR 8050), UPEM, UPEC, CNRS,
F-77454, Marne-la-Vall\'ee, France }\\
{\small $^2$ ENSAE-CREST-GENES
3, ave. P. Larousse
92245 MALAKOFF Cedex, FRANCE}\\
{\small $^3$ School of Mathematics and Statistics, Carleton University,
%1125 Colonel By Drive,
Ottawa, Ontario, K1S 5B6 Canada}
\date{}
}

\maketitle

\begin{abstract}
We consider the problem of  recovery of an unknown multivariate signal $f$ observed in a $d$-dimensional Gaussian white noise model of intensity $\e$.
We assume that $f$ belongs to a class of smooth functions ${\cal F}^d\subset L_2([0,1]^d)$ and has an  additive sparse structure
determined by the parameter $s$, the number of non-zero univariate components contributing to $f$.
We are interested in the case when $d=d_\e\to \infty$ as $\e\to 0$ and the parameter $s$ stays ``small'' relative to $d$.
With these assumptions, the recovery problem in hand becomes that of determining which sparse additive components are non-zero.
Attempting to reconstruct most non-zero components of $f$, but not all of them,
we arrive at the problem of almost full variable selection in high-dimensional regression.
For two different choices of ${\cal F}^d$,
we establish conditions under which almost full variable selection is possible, and
provide a procedure that gives almost full variable selection.
The procedure does the best (in the asymptotically minimax sense) in selecting most non-zero components of $f$.
Moreover, it is adaptive in the parameter $s$.
%Our results augment previous work in this area.

\medskip
\textbf{Keywords and phrases:} high-dimensional nonparametric regression; sparse additive
signals;
adaptive variable selection; exact and almost full selectors

\textbf{Mathematics Subject Classification:} 62G08, 62G20

\end{abstract}

%\noindent\textbf{AMS 2000 subject classifications:}

\section{Introduction}
In recent years, there has been much work on methods for variable selection in high dimensional settings;
refer, for example, to  \cite{Dalalyan, GJW2009, HHW, Ras} and references therein.
Among a variety of methods proposed, the lasso has become an important tool for sparse high-dimensional regression problems.
Motivated by the fact that finding the lasso solutions is computationally demanding,
Genovese et al. \cite{GJW2009} studied the relative statistical performance of the lasso and marginal regression,
which is also known as
simple thresholding,  for
sparse  high-dimensional regression problems. They found that marginal regression,
where each dependent variable is regressed separately
on each covariate, provides a good alternative to the lasso, and concluded that
their procedure merits further study.
Handling the problem of reconstruction in high dimensional regression,
Genovese et al. \cite{GJW2009} distinguished between the cases of exact, almost full, and no recovery.
Exact recovery refers to the situation where the set of all relevant components can be consistently recovered (asymptotically).
Almost full recovery stands for the possibility of having the number of misclassified components
negligibly small as compared to the
number of all relevant components.
The latter strategy requires milder restrictions on a statistical model and can be used in the situations where
exact recovery is impossible.  If neither exact nor almost full recovery can be achieved, we speak of `no  recovery' when the optimal risk is
as large as the number of relevant components and any recovery procedure fails completely.

Ingster and Stepanova \cite{IngStep2014}
extended the idea of Genovese et al. \cite{GJW2009} to the case of {nonparametric} regression.
Specifically, they addressed the problem of recovering sparse additive smooth signals observed in the continuous regression model
and showed that, asymptotically, as dimension increases indefinitely, exact variable selection is possible
and is provided by a suitable thresholding procedure. The procedure in \cite{IngStep2014} is optimal in the asymptotically minimax sense.
It is also free from the sparsity parameter and thus is adaptive.
At the same time, the more intricate problem of
almost full recovery in an adaptive setup remained unsolved. We shall treat this problem in the present paper.

Our setting is that of a multivariate signal $f\in{\cal F}^d\subset L_2([0,1]^d)=L_2^d$
corrupted by a Gaussian white noise of a given intensity $\e$:
\begin{gather}\label{model}
X_\e=f+\e W,
\end{gather}
where $W$ is a $d$-dimensional Gaussian white noise on $[0,1]^d$, $\e>0$ is a noise intensity, and ${\cal F}^d$
is a subset of $L_2^d$ that  consists of sufficiently smooth functions. In the present paper, two examples of ${\cal F}^d$
 will be considered.
In this model, the ``observation''
is the function $X_\e\colon L_2^d\to {\cal
G}$ taking its values in the set ${\cal G}$ of normal random
variables such that if $\xi=X_\e(\phi),\ \eta=X_\e(\psi),$ where
$\phi,\psi\in L_2^d$, then $\E(\xi)=(f,\phi),$ $\E(\eta)=(f,\psi),$
and $\Cov(\xi,\eta)=\e^2(\phi,\psi)$. For any $f\in L^d_2$, the
observation $X_\e$
determines the Gaussian measure $\Pb_{\e,f}$ on the Hilbert space $L_2^d$ with mean function $f$
and covariance operator $\e^2I $, where $I$ is the identity operator
(see \cite{IH.97, SK} for references). The expectation that corresponds to the probability measure $\Pb_{\e,f}$
is denoted by $\E_{\e,f}$.
In this paper,  the case of growing dimension
$d=d_\e\to \iy$ as $\e\to 0$ is studied.
It is well known  that the continuous model (\ref{model}) serves as a good approximation
to a more realistic equidistant sampling scheme with discrete Gaussian white noise.
In such an approximation, $\e^{-2}$ roughly corresponds to the number $n$ of observations per unit cube $[0,1]^d$.

An important problem in this context is to recover $f$ from noisy data.
Attempting to suppress the curse of dimensionality and complement the findings in \cite{IngStep2014}, we assume that  $f$ has an additive sparse structure.
Our goal is to study under what conditions and by means of what procedure
\textit{almost full} recovery of an additive sparse signal $f$ is possible.
In other words, we wish to correctly identify most non-zero %sparse additive
components of $f$.
In doing so, we aim at providing the procedure that, for the two function spaces ${\cal F}^d$ of our interest,
one consisting of functions of finite smoothness and the other consisting of functions of infinite smoothness, is optimal
in the asymptotically minimax sense.
In the almost full recovery regime, one can detect even smaller relevant components
but, unfortunately,  at the price of a loss in the rate.
Therefore constructing the corresponding procedure is technically more demanding
as compared to that in the exact recovery case.
To develop a good almost full recovery procedure, we will use results from
minimax hypothesis testing and minimax estimation theory.

To fix some notation  and assumptions,  let the signal $f$ in model (\ref{model}) be
of the form
(see, for example, \cite{GI} and \cite{IngStep2014})
\begin{gather*}%\label{adstruct}
f(\xb)=\sum_{j=1}^d \eta_j f_j(x_j),\quad \xb=(x_1,\ldots,x_d)\in[0,1]^d,\quad \eta=(\eta_1,\ldots,\eta_d)\in {\cal H}_{d,s},
\end{gather*}
where for a number $s\in\{1,\ldots,d\}$, called the \textit{sparsity parameter},
$${\cal H}_{d,s}=\{\eta=(\eta_1,\ldots,\eta_d): \eta_j\in\{0,1\}, 1\leq j\leq d,\, \sum_{j=1}^d \eta_j =s\}.$$
The $\eta_j$'s are non-random quantities taking values 0 and 1;
the case $\eta_j=1$ ($\eta_j=0$) corresponds to the situation when the component $f_j$ is active (non-active).
When $s=o(d)$ we speak of a \textit{sparse} additive signal $f$. In addition,
each component $f_j$ is assumed to be an element of a certain smooth function space ${\cal F}_\s\subset L_2[0,1]$ depending on a \textit{known}
parameter $\s>0$; two examples of ${\cal F}_{\s}$ under study are
introduced in Section 2.
Thus, the class of $s$-sparse multivariate signals of interest is
$${\cal F}^d_{s,\s}=\left\{f: f(\xb)=\sum_{j=1}^d \eta_j f_j(x_j),\;\int_0^1 f_j(x) \, dx=0, \;
f_j\in{\cal F}_{\s},\;1\leq j\leq d,\; \eta=(\eta_j)\in{\cal H}_{d,s}\right\},$$
where the components satisfy side condition that guarantees uniqueness,
and the signal recovery problem becomes that of determining which sparse additive components are non-zero.

In the context of variable selection, the problem of reconstruction of an additive function $f$ is now stated as follows.
For each component $f_j$ of a signal $f\in{\cal F}^d_{s,\s}$, consider testing the hypothesis of no signal $H_{0j}: f_j= 0$
versus the alternative $H_{1j}: f_j\in{\cal F}_\s(r_\e)$,
where for a positive family  $r_\e\to 0$
\begin{gather}\label{F}
{\cal F}_\s(r_\e)=\{g\in{\cal F}_\s: \|g\|_\s\leq 1,\;\|g\|_2\geq r_\e\},
\end{gather}
and $\|\cdot\|_\s$ is a norm on ${\cal F}_\s$. In this problem, a precise demarcation between
the signals that can be detected with error probabilities tending to 0 and the signals that
cannot be detected
is given in terms of a
\textit{detection boundary}, or \textit{separation rate},  $r_\e^*\to 0$ as $\e\to 0.$
For various function classes frequently used in minimax hypothesis testing,
sharp asymptotics for $r_\e^*$ are available (see, for example, \cite{I1993}).
The hypotheses $H_{0j}$ and $H_{1j}$
\textit{separate asymptotically} (that is,
the minimax error probability tends to zero)
if $r_\e/r_\e^*\to \iy$ as $\e\to 0$.  The hypotheses $H_{0j}$ and $H_{1j}$ \textit{merge asymptotically} (that is,
the minimax error probability tends to one)
if $r_\e/r_\e^*\to 0$ as $\e\to 0$.

When $H_{0j}$ and $H_{1j}$ separate asymptotically, we say that $f_j$ is \textit{detectable}.
If the hypotheses $H_{0j}$ and $H_{1j}$ separate (merge) asymptotically when
$\liminf r_\e/r_\e^*>1$ ($\limsup r_\e/r_\e^*<1$), the detection boundary $r_\e^*$ is said to be \textit{sharp}.
The knowledge of a sharp detection boundary $r_\e^*$ allows us to have a meaningful  problem
of testing $H_{0j}:f_j= 0$ versus $H_{1j}:f_j\in{\cal F}_\s(r_\e)$ by choosing $r_\e$ so that $\liminf_{\e\to 0}r_\e/r_\e^*>1$.
Otherwise, the function $f_j$ will be too ``small'' to be noticeable.

Let us agree to say that any measurable function $\eta^*=\eta^*(X_\e)$ taking values on $\{0,1\}^d$ is a \textit{selector}.
Following \cite{GJW2009} and \cite{IngStep2014}, we judge the quality of a selector $\eta^*$ of  vector $\eta\in {\cal H}_{d,s}$
by using the \textit{Hamming distance} on $\{0,1\}^d$, which counts the number of positions at which $\eta^*=(\eta_1^*,\ldots,\eta_d^*)$
and $\eta=(\eta_1,\ldots,\eta_d)$ differ:
$$|\eta^*-\eta|=\sum_{j=1}^d |\eta^*_j-\eta_j|.$$
%When dealing with the problem of a signal recovery, one
Following \cite{GJW2009}, we
distinguishe between  exact and almost full recovery. % (see, for example, \cite{GJW2009} and \cite{IngStep2014}).
Roughly, a selector $\eta^*=\eta^*(X_\e)$ is asymptotically \textit{exact} if its maximum risk is $o(1)$.
Likewise, a selector $\eta^*=\eta^*(X_\e)$ is asymptotically \textit{almost full} if its maximum risk is $o(s)$
with $s$ being the number of non-zero components $f_j$ of a  signal $f=\sum_{j=1}^d \eta_j f_j.$

Ingster and Stepanova \cite{IngStep2014} have obtained adaptive procedure that gives
asymptotically exact reconstruction of a $\s$-smooth signal $f\in{\mathcal F}^d_{s,\s}$ observed in a $d$-dimensional Gaussian white noise model.
A similar result for the space of infinitely-smooth functions is stated in this paper in Section 4.2 (see Theorems 1 and 2).
Although the selector in Section 4.2 is based on somewhat different statistics when compared to the one in  \cite{IngStep2014}, both selectors
have one common feature that their thresholds are free of the sparsity parameter $s$ and therefore automatically adapt themselves to its values.

The goal of this paper is three-fold. First, we find a sharp detection boundary %in terms of  $\sigma$,  $d$, and $s$
that allows us to separate detectable components of a signal $f\in{\mathcal F}^d_{s,\s} $ from non-detectable ones. Next, assuming that
all active components $f_j$
are detectable and that $s$ belongs to a set $\mathcal{S}_d$, which puts some mild restrictions on the range of $s$, we construct a selector $\eta^*=\eta^*(X_\e)$
with the property
\begin{gather}\label{er1}
\sup_{s \in \mathcal{S}_d} \sup_{\eta\in{\mathcal H}_{d,s}}\sup_{f\in{\mathcal F}^d_{s,\s}}s^{-1}\E_{f,\eta}|{\eta}^*-\eta|\to 0, \quad \mbox{ as } \e \to 0.
\end{gather}
Finally, we show that if at least one of the $f_j$'s is undetectable,
then
\begin{gather}\label{er2}
\liminf_{\e\to 0} \inf_{\ti{\eta}}\sup_{\eta\in{\mathcal H}_{d,s}}\sup_{f\in{\mathcal F}^d_{s,\s}}s^{-1}\E_{f,\eta}|\ti{\eta}-\eta|>0,
\end{gather}
that is, almost full recovery is impossible.
%By the symmetry of the model, %all the components $f_j$ are alike. Therefore
%if at least one of the $f_j$'s is not detectable,
%than all the components are not detectable.

The selector $\eta^*$ that satisfies (\ref{er1}) is said to provide asymptotically
\textit{almost full recovery} of a signal $f\in{\cal F}^d_{s,\s}$ in model (\ref{model});
its maximum risk is small relative to the number of non-zero components.
If, in addition, inequality (\ref{er2}) holds true,
then the selection procedure based on $\eta^*$ is the best possible (in the asymptotically minimax sense).
The notion of optimality that we use is borrowed from the minimax hypothesis testing theory.

In the present setup, adaptive (in $s$) variable selection in high dimensions presents several challenges.
First, one has to construct a good non-adaptive selector.
Second, having that selector available, one has to adapt it to unknown values of the parameter $s$.
It turns out that, when $s$ is known, both exact and almost full recovery can be achieved
by a suitably designed thresholding procedure (see Section 3.1 for details).
The problem of adaptation of this procedure to unknown values of $s$
was  tackled and solved in \cite{IngStep2014}, but in the case of exact recovery only.
Handling the same problem in the almost full recovery case will bring us in this paper to
the use of Lepski's method. This method  was proposed
for adaptive estimation in a Gaussian white noise model.
The reason why adaptive reconstruction of most relevant components of $f$ turns out to be more challenging
than adaptive reconstruction of all components of $f$
lies in the very nature of the thersholding procedure as defined in (\ref{ct}).
In contrast to the exact selector given by (\ref{etastar}) whose threshold is set regardless of the value of $s$,
thresholding in (\ref{ct}) does depend on $s$.

The paper is organized as follows.
In Section 2 we  present some general results of the asymptotically minimax hypothesis testing theory
and provide  details on their use for the two function spaces of our interest.
In Section 3 we translate the initial problem to the one in terms of the Fourier coefficients
and, for both function spaces in hand, obtain almost full selectors for a known sparsity parameter $s$.
In addition to that, we derive conditions under which almost full variable selection is possible.
Adaptive selectors for the function spaces in hand are developed in Section 4. To complete the picture,
we also introduce an adaptive  selection procedure that gives exact reconstruction for the space of analytic functions.
 Our main results
are stated in Section 4 and proved in Section 5.

\section{The building blocks}
As in \cite{IngStep2014}, the recovery problem under study will be connected to that of hypothesis testing.
Before stating and proving our main results, we shall discuss
some important tools of  minimax hypothesis testing
that will be used in the subsequence sections. For a complete exposition of the subject,
see \cite{IS.02} and the review papers \cite{I1993, I1993b, I1993c}.

\subsection{Extreme problem for ellipsoids: general case}
In asymptotically minimax hypothesis testing, when dealing with classes of smooth functions, the first common step is
to transform the initial problem involving a class of functions to the corresponding problem in the space
of Fourier coefficients. For this, let $\{\phi_k(x)\}_{k\in \mathbb{Z}}$ be the orthonormal basis in $L_2[0,1]$ given by
$$\phi_0(x)=1,\quad \phi_k(x)=\sqrt{2}\cos(2\pi kx),\quad \phi_{-k}(x)=\sqrt{2}\sin(2\pi kx),\quad k>0.$$
If $g\in L_2[0,1]$, then $g(x)=\sum_{k\in \mathbb{Z}}\t_{k} \phi_k(x)$, where
$\t_k=(g,\phi_k)$ is the $k$th Fourier coefficient of $g$, and $\|g\|_2^2= \sum_{k\in \mathbb{Z}}\t_k^2.$
Let  ${\cal F}_\s$ be a function space depending on a parameter $\s>0$ that is a subset of $L_2[0,1]$.
Suppose that $g\in{\cal F}_\s\subset L_2[0,1]$
is observed in a univariate Gaussian white noise of intensity $\e$, and we wish to test the null hypothesis
$H_0:g=0$ versus a sequence of alternatives $H_{1\e}: g\in  {\cal F}_\s(r_\e)$, where
the set ${\cal F}_\s(r_\e)$ is given by (\ref{F}).
For the two function spaces of interest, the norm of an element $g$ is expressed as
$\|g\|_\s^2= \sum_{k\in \mathbb{Z}}c_k^2\t^2_{k}$ with specified coefficients $c^2_k=c^2_k(\s)$
(see formulas (\ref{norma2}) and (\ref{norma1}) below).
In the sequence space of Fourier coefficients, the set  ${\cal F}_\s(r_\e)$ corresponds to
the ellipsoid in the space $l_2(\mathbb{Z})$ with semi-axis $c_k=c_k(\s)$ and a small neighbourhood of the point $\t=0$ removed:
\begin{gather}\label{Th}
\Theta_{\s}(r_\e)=\bigg\{\t=(\t_{k})_{k\in \mathbb{Z}}\in l_2(\mathbb{Z}):\sum_{k\in \mathbb{Z}}c_k^2\t^2_{k}\leq 1,\,
\sum_{k\in \mathbb{Z}}\t^2_{k}\geq r_\e^2\bigg\}.
\end{gather}

For constructing an asymptotically almost full selector, we shall need some facts from the minimax theory of hypothesis testing.
Denote by  $\t^*(r_\e)=(\t_k^*(r_\e))_{k\in \mathbb{Z}}$ the solution to the extreme problem
\begin{gather}\label{extr}
\frac{1}{2\e^4}\sum_{k\in\mathbb{Z}} \t_k^4 \to \inf_{\t\in \Theta_{\s}(r_\e)},
\end{gather}
and let $u_\e^2(r_\e)=u_\e^2(\Theta_\s(r_\e))$ be the value of the problem, that is,
$$u_\e^2(r_\e)=\frac{1}{2\e^4}\inf_{\t\in \Theta_{\s}(r_\e)}\sum_{k\in\mathbb{Z}}
\t_k^4=\frac{1}{2\e^4}\sum_{k\in\mathbb{Z}}\left( \t_k^*(r_\e) \right)^4.$$
The function $u_\e^2(r_\e)$ plays a key role in the minimax theory of hypothesis testing.
It controls the minimax total error probability and is used to set a cut-off point of the asymptotically minimax test procedure.
The detection boundary $r_\e^*$ in the problem of testing $H_0:\t=0$ versus $H_1:\t\in\Theta_\s(r_\e)$ is determined by
the relation $u_\e(r^*_\e)\asymp1.$
The function $u_\e(r_\e)$  is a non-decreasing function of the argument $r_\varepsilon$ which
possesses a kind of `continuity' property.
Namely, for any  $\epsilon>0$ there exist $\Delta>0 $ and $\varepsilon_0>0$ such that for any $\delta\in(0,\Delta)$
and $\varepsilon\in(0,\varepsilon_0)$,
\begin{gather}\label{ucont}
u_\e(r_\e)\leq u_\e((1+\delta)r_\e)\leq (1+\epsilon)u_\e(r_\e).
 \end{gather}
These and some other facts about $u_\e^2(r_\e)$ can be found in  \cite[Sec.\,3.2]{I1993} and \cite[Sec.\,5.2.3]{IS.02}).

For standard function spaces with the norm $\|g\|_\s$ defined (under the periodic constraints) in terms
of Fourier coefficients as $\|g\|_\s^2=\sum_{k\in\mathbb{Z}}\t_k^2 c_k^2$,
the form of the extremal sequence $(\t_k^*(r_\e)_{k\in\mathbb{Z}}$ in problem (\ref{extr}) as well as
the sharp asymptotics for $u_\e(r_\e)$ are available. Below we cite some relevant results for the two
function spaces ${\cal F}_\s$ of our interest: the Sobolev space of periodic $\s$-smooth function on $\mathbb{R}$ and
the space of periodic  functions on $\mathbb{R}$ that admit an analytic continuation to the strip around the real line.

\subsection{Extreme problem for  Sobolev ellipsoids}
Let ${\cal F}_\s$ with $\s>0$ denote the Sobolev space of $\s$-smooth $1$-periodic
functions on $\mathbb{R}$.
Define the norm $\|\cdot\|_\s$ on ${\cal F}_\s$  by the formula
\begin{gather}\label{norma2}
\|f\|_\s^2=\sum_{k\in \mathbb{Z}}\t_k^2 c_k^2,\quad c^2_k=c^2_k(\s)=(2\pi |k| )^{2\s},
\end{gather}
where $\t_k$ is the $k$th Fourier coefficient of $f$ with respect to  $\{\phi_k(x)\}_{k\in \mathbb{Z}}.$
If $\s$ is an integer, then under the periodic constraints (when the function admits 1-periodic $[\s]$-smooth
extension on the real line)  the norm as in (\ref{norma2}) corresponds to
\begin{gather*}
\|f\|_\s^2=\int_0^1\left(f^{(\s)}(x)\right)^2\, dx.
\end{gather*}

For a function $f\in{\cal F}_\s$ consider testing the hypothesis $H_{0}: f = 0$
versus the alternative $H_{1}: f\in{\cal F}_\s(r_\e)$,
where for a positive family  $r_\e\to 0$
\begin{gather*}
{\cal F}_\s(r_\e)=\{f\in{\cal F}_\s: \|f\|_\s\leq 1,\;\|f\|_2\geq r_\e\}.
\end{gather*}
Switching from Sobolev balls $\{f\in{\cal F}_\s: \|f\|_\s\leq 1\}$  to
 Sobolev ellipsoids
$\{\t\in l_2(\mathbb{Z}): \sum_{k\in\mathbb{Z}} c_k^2\t_k^2\leq 1\}$ %with
%$c^2_k=c^2_k(\s)=(2\pi |k| )^{2\s}$,
leads to the problem  of testing $H_0:\t=0$ versus $H_1:\t\in\Theta_\s(r_\e)$. The test procedure
that does the best in distinguishing between the latter two hypotheses is obtained by solving the
extreme problem (\ref{extr}) with
the semi-axes $c_k$  defined as in (\ref{norma2}); see Section 3 of \cite{I1993} for details.
The extremal sequence $(\t_k^*(r_\e)_{k\in\mathbb{Z}}$
satisfies
(see, for example,  \cite[\S\,3.2]{I1993} and Theorem 2 in \cite{IS.05}):
\begin{gather}\label{teta}
(\t_k^*(r_\e))^2\asymp r_\e^{2+1/\s} \left(1-\left({2\pi |k|}/{K_\e}\right)^{2\s}\right)   \;\; \mbox{for}\;\; 1\leq |k|\leq K_\e
\quad \mbox{and}\quad\t_k^*(r_\e)=0\;\;\mbox{otherwise},
\end{gather}
where
\begin{gather}\label{kes}
K_\e=\lfloor(4\s+1)^{1/(2\s)} r_\e^{-1/\s}\rfloor.
\end{gather}
The sharp asymptotics for $u_\e(r_\e)$ are of the form
(see \cite[\S\,4.3.2]{IS.02} and Theorems 2 and 4 in \cite{IS.05})
\begin{gather}\label{ta}
u_\e(r_\e)\sim C(\s) r_\e^{2+1/(2\s)} \e^{-2},\quad \e\to 0,
\end{gather}
where  (see, for example, p. 104 of \cite{I1993})
\begin{gather*}
C(\s)=2\s\left[  \left(1+\frac{1}{4\s}\right)\left( 1+4\s \right)^{1/(2\s)}\left(B\left(\frac{1}{2\s},2\right)\right)^{1/\s}  \right]^{-1},
\end{gather*}
and $B(\cdot,\cdot)$ is the Euler beta-function.

\subsection{Extreme problem for the ellipsoids of analytic functions} The following example of ${\cal F}_\s$ is also well known in nonparametric
estimation and hypothesis testing.
Let ${\cal F}_\s$ with $\s>0$ be the class of 1-periodic functions
$f$ on $\mathbb{R}$ admitting a continuation to the strip $S_\s=\{z=x+iy: |y|\leq \s\}\subset\mathbb{C}$
such that $f(x+i y)$ is analytic on the interior of $S_\s$, bounded on $S_\s$ and
$$\int_{0}^1|f(x\pm i\s)|^2\, dx<\infty.$$ %
Let the norm $ \|\cdot \|_{1,\s}$ on ${\cal F}_\s$ be given by (see, for example, \cite{GL1996})
$$\|f\|^2_{1,\s}=\int_0^1\left( {\rm Re} f(x+i\s)  \right)^2 \, dx. $$%,\quad f\in {\cal F}_\s.$$
 In terms of the Fourier coefficients, the squared norm $\|f\|^2_{1,\s}$ takes the form
\begin{gather*}
\|f\|_{1,\s}^2=\sum_{k\in \mathbb{Z}}\t_k^2 c_k^2,\quad c^2_k=c^2_k(\s)=\cosh^2(2\pi \s k ).
\end{gather*}
In view of the relations
$$\exp(|x|)\leq 2\cosh(x)\leq 2\exp(|x|),\quad x\in \mathbb{R},$$ we may also consider an equivalent norm $\|\cdot\|_{\s}$ defined as
\begin{gather}\label{norma1}
\|f\|^2_{\s}= \sum_{k\in \mathbb{Z}}\t_k^2 c_k^2,\quad c_k=c_k(\s)=\exp(2\pi \s |k|).
\end{gather}
We have chosen to deal with the latter norm as it is easier to study.

The ball $\{f\in{\cal F}_\s: \|f\|_\s\leq 1\}$ corresponds  to
the ellipsoid
$\{\t\in l_2(\mathbb{Z}): \sum_{k\in\mathbb{Z}} c_k^2\t_k^2\leq 1\}$
with the semi-axes $c_k$  defined as in (\ref{norma1}).
Thus translating the problem of testing $H_0:f=0$ versus $H_1: f\in{\cal F}_\s(r_\e)$
to the one in terms of Fourier coefficients brings us to testing $H_0:\t=0$ versus $H_1:\t\in\Theta_\s(r_\e)$. The
asymptotically minimax test procedure
that  distinguishes between these two hypotheses is obtained by solving the
extreme problem (\ref{extr}) with
the semi-axes $c_k$  defined as in (\ref{norma1}).
The elements  of the extremal sequence $(\t_k^*(r_\e)_{k\in\mathbb{Z}}$
in problem (\ref{extr}) with the semi-axis $c_k$ as above %defined in (\ref{norma1})
may be taken as constants (independent of $k$) satisfying as $\e\to 0$  (see, for example, Section 3 in \cite{I1993})
\begin{gather}\label{elements}
\t_k^*(r_\e)\asymp r_\e\log^{-1/2}(r_\e^{-1})\;\;\mbox{for}\; \;1\leq |k|\leq K_\e\quad \mbox{and}\quad
\t_k^*(r_\e)=0\;\;\mbox{otherwise},
\end{gather}
where
\begin{gather}\label{kes1}
K_\e=\lfloor(2\pi \s)^{-1}\log (r_\e^{-1})\rfloor,
\end{gather}
and we have
\begin{gather}\label{ta1}
u_\e(r_\e)\sim \left(\frac{r_\e}{\e}  \right)^2\frac{(2\pi \s)^{1/2}}{\log^{1/2} (r_{\e}^{-1})}.
\end{gather}

Formulas (\ref{elements})--(\ref{ta1}), as well as formulas (\ref{teta})--(\ref{ta}),
will be employed to construct almost full selectors for the two function spaces under study.

\section{Variable selection in a sequence space model}
By sufficiency,  the problem of recovering $f$ observed in the Gaussian white noise model can be transformed to an
equivalent problem in a sequence space model.
Acting as in \cite{IngStep2014}, for the index $l\in \mathbb{Z}^d$
whose $j$th component is equal to $k$ and the other components are equal to zero,
define the function
$$\phi_{j,k}(\xb)=\phi_l(\xb)=\phi_k(x_j),\quad\xb=(x_1,\ldots,x_d)\in[0,1]^d,\quad 1\leq j\leq d,\quad k\in \mathbb{Z},$$
and denote by $\t_{j,k}=(f,\phi_{j,k})=\int_0^1\phi_k(x) f_j(x)\, dx$ the $k$th
Fourier coefficient  of the $j$th component  $f_j$ of a signal $f=\sum_{j=1}^d \eta_j f_j$.
Consider the sequence space model
\begin{gather}\label{model1}
X_{j,k}=\eta_j\t_{j,k}+\e \xi_{j,k}, \quad \xi_{j,k}\stackrel{i.i.d.}\sim{{N}}(0,1),\quad 1\leq j\leq d, \quad k\in \mathbb{Z},
\end{gather}
where $X_{j,k}=X_\e(\phi_{j,k})$ are the empirical Fourier coefficients and the collection
$(\eta_1\t_1,\ldots,\eta_d\t_d)$ consists of sequences $\eta_j\t_j=(\eta_j\t_{j,k})_{k\in \mathbb{Z}}$
such that
$(\eta_j)\in {\mathcal H}_{d,s}$ and for all $1\leq j\leq d$,
\begin{gather}\label{set1}
\t_j=(\t_{j,k})\in l_2(\mathbb{Z}),\quad\sum_{k\in \mathbb{Z}}c_k^2\t^2_{j,k}\leq 1.
\end{gather}
In this paper we have chosen to deal with the latter model,
which is technically more convenient.
Although the set of $\t_j$s in (\ref{set1}) involves an orthogonal system in $L_2^d$, the results on minimax errors and risks
do not depend on the choice of this orthogonal system because
the random variables $X_{j,k}$,
which generate a sufficient $\s$-algebra for $f\in{\cal F}^d_{s,\s}$, are independent normal $N(\eta_j\t_{j,k},\e^2)$.
Thus the distribution of $\{X_{j,k}\}$ depends on the Fourier coefficients
$\t_{j,k}$ of $f$ with respect to the system $\{\phi_{j,k}\}$ but not on the choice of $\{\phi_{j,k}\}$.
Using a suitable finite collection of the random variables  $X_{j,k}$ as defined in (\ref{model1}), we wish
to construct an optimal selection procedure
that recovers most non-zero components of $(\eta_1\t_1,\ldots,\eta_d\t_d)$, but not all of them.

\subsection{Almost full variable selection in the non-adaptive case}

We first consider a non-adaptive setup when the sparsity parameter $s$ is known.
When dealing with the problem of variable selection
in model (\ref{model1}), we make use of
the statistics, cf. asymptotically minimax test statistics in Section 3.1 of \cite{I1993},
\begin{gather}\label{tj}
t_j=t_j(s)=\sum_{1\leq |k|\leq K_\e}\omega_k(r_{\e}^*(s))\left[\left( \frac{X_{j,k}}{\e} \right)^2-1  \right],\quad j=1,\ldots,d,
\end{gather}
where for any $r_\e>0$ the weight functions  $\omega_k(r_\e)$ are given by the formula
$$
\omega_k(r_\e)=\frac{1}{2\e^2}\frac{(\t_k^*(r_\e))^2 }{ u_\e(r_\e)},\quad 1\leq |k|\leq K_\e,
$$
and the number $r_{\e}^*(s)>0$ is the solution of the equations
\begin{gather}\label{re}
\frac{u_\e(r^*_{\e}(s))}{\sqrt{2\log (d/s)}}=1.
\end{gather}
For both function spaces of interest, the quantities $K_\e$, $\t_k^*(r_\e)$, and $ u_\e(r_{\e})$  in formula (\ref{tj}) are specified in Section 2.
The sparsity parameter
$s\in\{1,2,\ldots,d\}$ is assumed to be small relative to $d$, that is, $s=o(d)$.
Note that the weights $\omega_k(r_\e)$ are normalized to have $\sum\limits_{1\leq |k|\leq K_\e}\omega_k^2(r_\e)=1/2$.

Now we define a non-adaptive \textit{almost full  selector} to be
\begin{gather}\label{ct}
\check{\eta}=(\check{\eta}_1,\ldots,\check{\eta}_d),\quad \check{\eta}_j=\mathbb{I}\left( t_j>\sqrt{2\log (d/s)+\delta\log d}\right),\quad j=1,\ldots,d,
\end{gather}
where  $\delta=\delta_\e>0$ satisfies
\begin{gather}\label{deltacond0}
\delta\to 0\quad\mbox{and}\quad \delta\log d\to \infty,\quad \mbox{as}\; \;\e\to 0.
\end{gather}
The arguments as in the proof of Theorem 1 show that for Sobolev ellipsoids,
under the conditions, cf. (\ref{dr}),
$$\log d=o(\e^{-2/(2\s+1)}),\quad\liminf_{\e\to 0}\frac{u_\e(r_\e)  }{\sqrt{2\log (d/s)}}>1,$$
the selector $\check{\eta}$ reconstructs almost all relevant components of a vector $\eta\in {\mathcal H}_{d,s}$, and hence
asymptotically provides almost full recovery of a signal $f\in{\mathcal F}^{d}_{s,\s}$ in model (\ref{model}).

To illustrate the difference between exact and almost full reconstruction in adaptive settings,
assume that ${\cal F}_\s$ is the Sobolev space.
In this case,
a selector (see Section 3.1 of \cite{IngStep2014} with $s$ in place of $d^{1-\beta}$)
\begin{gather}\label{etastar}
{\eta}^*=({\eta}^*_1,\ldots,{\eta}^*_d),\quad {\eta}^*_j=\mathbb{I}\left( t^*_j>\sqrt{(2+\delta)\log d}\right),\quad j=1,\ldots,d,
\end{gather}
where the statistics $t_j^*$ are defined similar to the $t_j$ as in (\ref{tj}) with the relation
\begin{gather*}
\frac{u_\e(r^*_{\e}(s))}{\sqrt{2\log d}+\sqrt{2\log s }}=1
\end{gather*} instead of (\ref{re}),
turns our to be a non-adaptive \textit{exact selector}, as long as
\begin{gather}\label{dr}
\log d=o(\e^{-2/(2\s+1)})\quad\mbox{and}\quad \liminf_{\e\to 0}\frac{u_\e(r_\e)  }{\sqrt{2\log d}+\sqrt{2\log s}}>1.
\end{gather}
Under the above conditions, the procedure based on $\eta^*$ selects correctly all non-zero components of a vector $\eta\in {\mathcal H}_{d,s}$, and hence
provides exact recovery of a signal $f\in{\mathcal F}^d_{s,\s}$ in model (\ref{model}).

Contrasting with formula (\ref{etastar}), the threshold
in (\ref{ct}) is set at a lower level and is dependent on the parameter $s$. The latter fact
makes the idea of adaption suggested in \cite{IngStep2014} for the
exact reconstruction case invalid (see Section 3.3. for details). In the next section, we obtain the desired adaptive %almost full
selector by using Lepski's method. Before doing that, we provide conditions on
$d$ as a function of $\e$ under which the thresholding procedure (\ref{ct}), as well as its adaptive version
introduced in Section 4.1, gives asymptotically almost full reconstruction
of a function $f\in{\mathcal F}^d_{s,\s}$.

\subsection{Conditions for almost full variable selection}
Consider now the question of determining conditions on  $d$ as a function of $\e$  under which almost full variable selection is possible.
Violation of these conditions will lead to entirely different selection strategies.

In the sequence space of Fourier coefficients, consider testing the null hypothesis $H_{0j}: \t_j=0$ versus the alternative  $H_{1j}:\t_j\in \Theta_\s(r_{\e})$,
where the set $\Theta_\s(r_{\e})$ is given by (\ref{Th}). It is easy to see that
under the null hypothesis $H_{0j}$, we have (see, for example, Section 4.1 of \cite{IngStep2014})
 $$\E_0(t_{j})=0,\quad \Var_0(t_{j})=1,$$
while under the alternative $H_{1j}$,
where for all sufficiently small $\e$ a small parameter $r_\e>0$ satisfies $r_\e/r_{\e}^*(s)>1$,
\begin{eqnarray}
\E_{\t_j}(t_{j})&=&\e^{-2}\sum_{1\leq |k|\leq K_\e}\omega_k(r^*_{\e}(s))\t_{j,k}^2\geq u_\e(r^*_{\e}(s)),\label{mo}\\
\Var_{\t_j}(t_{j})&=&1+O(\E_{\t_j}(t_{j})\max_{1\leq |k|\leq K_\e}\omega_k(r^*_{\e}(s))). \nonumber %\label{disp}
\end{eqnarray}
Furthermore, under the above restrictions on $r_\e$ and $r_{\e}^*(s)$ the following result holds (in case of Sobolev spaces, see Proposition 7.1 in \cite{GI} and
Lemma 1 in \cite{IngStep2014}; in case of the space ${\cal F}_\s$ of analytic functions, the proof is similar to that of Sobolev spaces).

Let the quantity $T=T_\e\to -\iy$ and the weight functions $\omega_k(r^*_{\e}(s))$ as in {\rm (\ref{tj})} be such that
as $\e\to 0$
\begin{gather}\label{usl}
 T\max_{1\leq |k|\leq K_\e}\omega_k(r^*_{\e}(s))\to 0\quad\mbox{and}\quad \E_{\t_j}( t_{j}(s)) \max_{1\leq |k|\leq K_\e}\omega_k(r^*_{\e}(s)) \to 0.
 \end{gather}
Then as $\e\to 0$
\begin{gather}\label{expon1}
\Pb_{0}(t_{j}\leq T)\leq \exp\left(-\frac{T^2}{2} (1+o(1))\right),
\end{gather}
and for all $j=1,\ldots, d$, uniformly in $\t_j\in \Theta_\s(r_{\e})$,
\begin{gather}\label{expon2}
\Pb_{\t_j}(t_{j}-\E_{\t_j}( t_{j})\leq T)\leq \exp\left(-\frac{T^2}{2} (1+o(1))\right).
\end{gather}

For both function spaces ${\cal F}_\s$ of our interest, the exponential bounds (\ref{expon1}) and (\ref{expon2})
will be applied below to the quantity $T=T_\e\to-\infty$ of order $O(\log^{1/2}d)$.
This observation and assumption (\ref{re}) transform requirement (\ref{usl}) into
\begin{gather}\label{cond}
{\log^{1/2} d}\max_{1\leq |k|\leq K_\e}\omega_k(r_\e^{*}(s))\to 0,\quad \e\to 0,
\end{gather}
Condition (\ref{cond}) gives a restriction on the growth of $d=d_\e$ ensuring that
the selection procedure works as designed.
Indeed, as shown in Section 4.1 in \cite{IngStep2014}, for the Sobolev space of $\s$-smooth functions, one has
\begin{gather*}%\label{omega}
\omega_k(r_\e) \asymp r_\e^{1/(2\s)}\;\; \mbox{for }\;\; 1\leq |k|\leq K_\e,
\end{gather*}
and
\begin{gather*}
r_\e^{*}(s)\asymp \left(\e \log d\right)^{\s/(4\s+1)}.
\end{gather*}
Therefore condition  (\ref{cond}) is fulfilled when
\begin{gather}\label{c1}
\log d=o(\e^{-2/(2\s+1)})
\end{gather}

In case of the space ${\cal F}_\s$ of analytic functions, one has
\begin{gather}\label{omegaanalytic}
\omega_k(r_\e) \asymp \log^{-1/2}(r_\e^{-1})\;\; \mbox{for }\;\; 1\leq  |k|\leq K_\e,
\end{gather}
and, in view of (\ref{ta}) and (\ref{re}), the quantity $r_{\e}^*(s)$ satisfies
$$\log^{1/2}d\asymp \left(\frac{r_\e^{*}(s)}{\e}  \right)^2\log^{-1/2} \left((r_{\e}^{*}(s))^{-1}\right),$$
implying
\begin{gather*}
r_\e^{*}(s)\asymp \e \log ^{1/4} (d)\log ^{1/4}((r_\e^{*}(s))^{-1}).
\end{gather*}
Therefore $\log \left(\left(r_\e^{*}(s)\right)^{-1}\right)\sim \log (\e^{-1})$, and (see (\ref{omegaanalytic}))
$$\omega_k(r_\e^{*}(s))\asymp\log^{-1/2} (\e^{-1}).$$
From this, the technical condition (\ref{cond})
holds true when, cf. formula (\ref{omegaanalytic}),
\begin{gather}\label{condition}
\log d=o(\log(\e^{-1})),\quad \e\to 0.
\end{gather}

\section{Main results}
In this section, we consider a more realistic problem when the sparsity parameter $s$ is \textit{unknown}.
We derive conditions under which almost full variable selection is possible, and
construct a selector for which the Hamming distance is much smaller than the number of relevant components (see Theorems 3 and 5).
Our selector is adaptive in the sparsity parameter $s$ and is unimprovable in the asymptotically minimax sense (see Theorems 4 and 6).
In addition to that, in Section 4.2  we provide asymptotically exact selection procedure for the space of analytic
functions  that is adaptive in the sparsity parameter $s$.

\subsection{Almost full variable selection in the adaptive case}
In this subsection, the selector $\check{\eta}$  as in (\ref{ct}) will be used to obtain the corresponding adaptive procedure.
To avoid losses due to adaptation, we will have to limit the range of the possible values of $s$.
%in order to have an adaptive procedure
%Typically, adaptation incur some (not big) losses.
%The price  paid for adaptation here reveals itself
%in limiting the range of the possible values of $s$.
Namely, we assume that for some constants $0<c<C<1$
\begin{gather}\label{ogran}
c\leq\liminf_{d\to \infty}(\log s/\log d)\leq \limsup_{d\to \infty}(\log s/\log d)\leq C,
\end{gather}
and define the set
$$\mathcal{S}_d=\{s\in\{1,\ldots,d\}\;\mbox{is such that condition}\; (\ref{ogran})\; \mbox{holds}\}$$
over which the adaptive selector that we propose yields  almost full selection.
The restriction on $s$ as in (\ref{ogran}) is relatively mild. For instance, any
$s=d^{1-\beta}$ with $\beta\in[b, B]$ for some constants $0<b<B<1$ belongs to $\mathcal{S}_d$.

To construct the desired selector, for some $\Delta=\Delta_d>0$ and $M=\lceil (C-c)/\Delta\rceil+1$, pick grid points over the interval $(1,d)$:
\begin{gather}\label{nodes}
s_1=d^c,\quad s_m=s_{m-1}d^{\Delta}=s_1d^{(m-1)\Delta },\quad  2\leq m\leq M,
\end{gather}
and assume that
\begin{gather}\label{deltacond}
\Delta\to 0,\quad \Delta \log d\to 0,\quad\mbox{as}\;\; d\to \infty,
\end{gather}
yielding  $d^{\Delta}\leq \rm{const}$ for all large enough $d$.
For each $m=1,\ldots, M$, let the parameter $r_{\e}^*(s_m)>0$ be determined by the equation, cf. (\ref{re}),
 \begin{gather*}%\label{rstar}
\frac{u_\e(r_{\e}^*(s_m))}{\sqrt{2\log (d/s_m)}}=1,
\end{gather*}
where, depending on a type of
the ellipsoid $\Theta_\s(r_\e)$ we are dealing with, the function $u_\e(r_\e)$ satisfies either (\ref{ta})  or (\ref{ta1}).

Similar to the case of known $s$, consider weighted chi-square type statistics, cf. (\ref{tj}),
\begin{gather*}%\label{tjm}
t_{j}(s_m)=\sum_{1\leq |k|\leq K_\e}\omega_k(r_{\e}^*(s_m))\bigg[\Big( \frac{X_{j,k}}{\e} \Big)^2-1  \bigg],\quad j=1,\ldots,d,\quad m=1,\ldots,M.
\end{gather*}
with weight functions
\begin{gather*}%\label{wk}
\omega_k(r_{\e}^*(s_m))=\frac{1}{2\e^2}\frac{(\t_k^*(r_{\e}^*(s_m))^2 }{ u_\e(r_{\e}^*(s_m))},\quad 1\leq |k|\leq K_\e,
\end{gather*}
possessing the property $\quad \sum_{1\leq |k|\leq K_\e}\omega^2_k(r_{\e}^*(s_m))=1/2$.
The values of $\t_k^*(r_{\e}^*(s))$ and $K_\e$ depend on the function space under consideration.
For the Sobolev space in hand, $\t_k^*(r_{\e}^*(s))$ and $K_\e$ are as in (\ref{teta}) and (\ref{kes});
for the space of analytic functions,  $\t_k^*(r_{\e}^*(s))$ and $K_\e$ are as in (\ref{elements}) and (\ref{kes1}).

Next, for all $j=1,\ldots,d$ and $m=1,\ldots,M$, set
$${\wh{\eta}}_{j}(s_m)=\mathbb{I}\left(t_{j}(s_m)>\sqrt{2\log (d/s_m)+\delta\log d}   \right),$$
where $\delta=\delta_\e>0$ satisfies (\ref{deltacond0}),
and define an adaptive selector of a vector $\eta\in{\mathcal H}_{d,s}$ by the formula
\begin{gather}\label{etah}
\wh{\eta}(s_{\wh{m}})=(\wh{\eta}_1(s_{\wh{m}}),\ldots,\wh{\eta}_d(s_{\wh{m}})),
\end{gather}
where
$\wh{m}$ is chosen by Lepski's method (see Section 2 of \cite{Lepski1990}) as follows:
\begin{gather*}
\wh{m}=\min\left\{  1\leq m\leq M\,:\,|\wh{\eta}(s_m)-\wh{\eta}(s_i)|\leq v_{i}\;\;\mbox{for all}\;\; i\geq m   \right\}.
\end{gather*}
Here the quantities $v_i=v_{i,d}$ are set to be
$$v_i={s_{i}}/{\tau_d},\quad m\leq i\leq M,$$
with a sequence of numbers $\tau_d\to \infty$ satisfying (recall that $d=d_\e\to \infty$ as $\e\to 0$)
\begin{gather*}
{\tau_d}=o\left({\min(\log d, d^{\delta/2})}\right),\quad \mbox{as}\;\; \e\to 0.
\end{gather*}

Algorithmically, Lepski's procedure for choosing $\hat{m}$ works as follows.
We start by setting $\hat{m}=M$ and attempt to decrease the value of $\hat{m}$ from $M$
to $M-1.$ If $|\wh{\eta}(s_{M-1})-\wh{\eta}(s_M)|\leq v_{M}$, we set $\hat{m}=M-1$;
otherwise, we keep $\hat{m}$ equal to $M$. In case $\hat{m}$ is decreased to $M-1$, we continue the process attempting to
decrease it further. If $|\wh{\eta}(s_{M-2})-\wh{\eta}(s_{M-1})|\leq v_{M-1}$ and
$|\wh{\eta}(s_{M-2})-\wh{\eta}(s_M)|\leq v_{M}$,  we set $\hat{m}=M-2$; otherwise, we keep $\hat{m}$ equal to $M-1$; and so on.
Notice that by construction $v_M\geq v_{M-1}\geq \ldots \geq v_1.$

\subsection{Exact variable selection for analytic functions in the adaptive case}
The problem of adaptive reconstruction of sparse additive functions in the Gaussian white noise model
was studied in the only case of $\s$-smooth functions, see \cite{IngStep2014}.
 Before handling the problem of almost full variable selection in adaptive settings,
we complement the findings in \cite{IngStep2014} by presenting an adaptive exact selector
for the space of analytic functions. The strategy is similar to the one suggested  in \cite{IngStep2014} for $\s$-smooth functions, but the parameters of the statistics and the condition on the dimension $d$ are different.

Consider a sequence space model that corresponds to the Gaussian white noise model with $f$ from the class of analytic functions ${\cal F}_\s$  as defined in Section 2.3.
Let $1<s_1<s_1<\ldots <s_M<d$ be the grid of points as in (\ref{nodes}).
For any $m=1,\ldots, M$, let the parameter $r_{\e,m}^*>0$ be determined by the equation
 \begin{gather*}%\label{rstar}
\frac{u_\e(r^*_{\e,m})}{\sqrt{2\log d}+\sqrt{2\log s_m}}=1.
\end{gather*}
Consider weighted chi-square type statistics
\begin{gather*}%\label{tjm}
t_{j,m}=\sum_{1\leq |k|\leq K_\e}\omega_k(r_{\e,m}^*)\bigg[\Big( \frac{X_{j,k}}{\e} \Big)^2-1  \bigg],\quad j=1,\ldots,d,\quad m=1,\ldots,M.
\end{gather*}
with weight functions
\begin{gather*}%\label{wk}
\omega_k(r^*_{\e,m})=\frac{1}{2\e^2}\frac{(\t_k^*(r^*_{\e,m}))^2 }{ u_\e(r^*_{\e,m})}
\end{gather*}
obeying the normalization condition $\sum_{k\in \mathbb{Z}}\omega^2_k(r^*_{\e,m})=1/2.$
Next, for all $j=1,\ldots,d$ and $m=1,\ldots,M$, set
$$\eta_{j,m}=\mathbb{I}\left(t_{j,m}>\sqrt{(2+\delta)(\log d+\log M)}   \right),$$
and define an \textit{adaptive exact selector} $\eta^{**}$ of a vector $\eta\in{\mathcal H}_{d,s}$ by the formula (see formula (18) in \cite{IngStep2014})
\begin{gather}\label{etahstar}
\eta^{**}=(\eta_1^{**},\ldots,\eta_d^{**}),\quad \eta_j^{**}=\max_{1\leq m\leq M}\eta_{j,m},\quad j=1,\ldots,d.
\end{gather}
The idea behind the selector $\eta^{**}$ is as follows.
The $j$th component of a signal is viewed active if at least one of the statistics $t_{j,m}$, $m=1,\ldots,M$, detects it.
Therefore, thinking of $\eta_{j,m}$ and $\eta_j^{**}$ as test functions, we get that the probability of having $\t_j$ incorrectly undetected does not exceed the
respective probability with the $\eta_{j,m}$ test, where $s_m$ is
 close to the true (but unknown) value of $s$. Furthermore, the probability that $\eta_j^{**}$ incorrectly detects $\t_j$ is less than the sum of the respective probabilities
for the $\eta_{j,m}$ tests over all $m=1,\ldots,M,$
and is small by the choice of threshold.

Let the set $\Theta_{\s,d}(r_\e)$ be as in ({\ref{thd}) with the coefficients $c_k$ given by (\ref{norma1}).
The following two theorems,  whose proofs are similar to those of Theorems 3 and 4 in \cite{IngStep2014}, hold true.

\bigskip
\noindent \textbf{Theorem 1.} \textit{Let $s\in\{1,\ldots,d\}$ be such that $s=o(d)$.
Assume that $\log d=o(\log\e^{-1})$ and that the quantity
$r_\e=r_\e(s)>0$ satisfies
$$\liminf_{\e\to 0}\frac{u_\e(r_\e)  }{\sqrt{2\log d}+\sqrt{2\log s}}>1.$$
Then  as $\e\to 0$
$$\sup_{\eta\in{\mathcal H}_{d,s}}\sup_{\t\in\Theta_{\s,d}(r_\e)}\E_{\eta,\t}|\eta-{\eta}^{**}|\to 0,$$
where ${\eta}^{**}$ is the selector of vector $\eta $ as defined in {\rm (\ref{etahstar})}.}

\bigskip

\noindent
\textbf{Theorem 2.} \textit{Let $s\in\{1,\ldots,d\}$ be such that $s=o(d)$.
Assume that $\log d=o(\log\e^{-1})$ and that the quantity $r_\e=r_\e(s)>0$ satisfies
$$\limsup_{\e\to 0}\frac{u_\e(r_\e)  }{\sqrt{2\log d+\sqrt{2\log s}}}<1.$$
Then
$$\liminf_{\e\to 0}\inf_{\ti{\eta}}\sup_{\eta\in{\mathcal H}_{d,s}}\sup_{\t\in\Theta_{\s,d}(r_\e)} \E_{\eta,\t}|\eta-\ti{\eta}|> 0,$$
where the infimum is over all selectors $\ti{\eta}$ of vector $\eta$ in model {\rm (\ref{model1})}.}

\bigskip
\noindent \textbf{Remark 1.}
The sharp detection boundary in Theorems 1 and 2 which makes it possible to decide on whether we are in a position to
proceed further with  variable selection or not,
is determined in terms of the function $u_\e(r_\e)$ with sharp asymptotics as in  (\ref{ta}) and (\ref{ta1}).
The use of $u_\e(r_\e)$ instead of $r_\e$  makes it easy to build a bridge between variable selection in  Gaussian white noise setting
and variable selection in regression setting as studied in
Sec. 4 of \cite{GJW2009}. In addition, using $u_\e(r_\e)$ instead of $r_\e$ makes the statement of detectability condition precise.
By `continuity' of $u_\e(r_\e)$ as cited in (\ref{ucont}), the conditions of Theorems 1 and 2
that separate detectible components from undetectable ones can be written in a usual form
$\liminf_{\e\to 0}r_\e/r_\e^*>1$ and $\limsup_{\e\to 0}r_\e/r_\e^*<1$, where for Sobolev ellipsoids
the sharp detection boundary $r_\e^*$ is found explicitly from (\ref{ta}),
and for the ellipsoids of analytic functions it %$r^*_\e$
is found implicitly from (\ref{ta1}).
Similar remark applies to Theorems 3 to 6 stated in Section 4.3 and 4.4,
%The reason for using $u_\e(r_\e)$ instead of $r_\e$ is to make the statement of detectability condition precise.

\subsection{Almost full variable selection for Sobolev balls} Consider the set $\Theta_{\s}(r_\e)$ as in (\ref{Th})
with the coefficients $c_k$ given by (\ref{norma2}), and define the set
\begin{multline}\label{thd}
\Theta_{\s,d}(r_\e)
=\bigg\{\t=(\t_j): \t_j=(\t_{j,k})\in l_2(\mathbb{Z}),\
 \sum_{k\in\mathbb{Z}}c_k^2\t_{j,k}^2\leq 1, \
\sum_{k\in\mathbb{Z}} \t_{j,k}^2\geq r_\e^2,\
 1\leq j\leq d \bigg \}.%\label{ell11}
 \end{multline}
Let $\wh{\eta}(s_{\wh{m}})$ be the selector given by {\rm (\ref{etah})} based on the
statistics $t_j(s_m)$ as in (\ref{tj}), where the quantities  $\t_k^*(r_\e) $, $K_\e$, and $u_\e(r_\e)$
are specified by formulas (\ref{teta}), (\ref{kes}), and (\ref{ta}), respectively.
The following  theorem holds.

\bigskip
\noindent \textbf{Theorem 3.} \textit{Let  $s\in\{1,\ldots,d\}$ be such that {\rm (\ref{ogran})} holds true.
Assume  that $\log d=o(\e^{-2/(2\s+1)})$ and that the quantity
$r_\e=r_\e(s)>0$ satisfies
$$\liminf_{\e\to 0}\frac{u_\e(r_\e)  }{\sqrt{2\log (d/s)}}>1.$$
Then as $\e\to 0$}
$$\sup_{s\in\mathcal{S}_d}\sup_{\eta\in{\mathcal H}_{d,s}}\sup_{\t\in\Theta_{\s,d}(r_\e)}s^{-1}\E_{\eta,\t}|\wh{\eta}(s_{\wh{m}})-\eta|\to 0.$$

\bigskip

Theorem 3 says that if all the hypotheses $H_{0j}: \t_j\equiv0 $ and $H_{1j}:\t_j\in\Theta_\s(r_\e)$, $j=1,\ldots,d$,
separate asymptotically, then the selection procedure based on $\wh{\eta}(s_{\wh{m}})$
reconstructs almost all non-zero components of a vector $\eta\in {\mathcal H}_{d,s}$, and thus
provides almost full recovery of $(\eta_1\t_1,\ldots,\eta_d \t_d)$, uniformly in $\mathcal{S}_d$,
${\mathcal H}_{d,s}$, and $\Theta_{\s,d}(r_\e)$.

The next result shows that if the detectability condition is not met,
almost full selection is impossible.

\bigskip
\noindent
\textbf{Theorem 4.} \textit{Let  $s\in\{1,\ldots,d\}$ be such that $s=o(d)$.
Assume that $\log d=o(\e^{-2/(2\s+1)})$ and that the quantity $r_\e=r_\e(s)>0$ satisfies
$$\limsup_{\e\to 0}\frac{u_\e(r_\e)  }{\sqrt{2\log (d/s)}}<1.$$
Then
$$\liminf_{\e\to 0}\inf_{\ti{\eta}}\sup_{\eta\in{\mathcal H}_{d,s}}\sup_{\t\in \Theta_{\s,d}(r_\e)}s^{-1} \E_{\eta,\t}|\eta-\ti{\eta}|> 0,$$
where the infimum is over all selectors $\ti{\eta}$ of a vector $\eta$ in model {\rm (\ref{model1})}.}

\subsection{Almost full variable selection for analytic functions} The results similar to Theorems 3 and 4 hold true for the space of analytic functions.
Namely, consider the sets $\Theta_{\s}(r_\e)$ and $\Theta_{\s,d}(r_\e)$ as in (\ref{Th}) and (\ref{thd})
with the coefficients $c_k$ given by (\ref{norma1}).
Again, let $\wh{\eta}(s_{\wh{m}})$ be the selector defined by {\rm (\ref{etah})} based on the
statistics $t_j(s_m)$ as in (\ref{tj}), but the quantities  $\t_k^*(r_\e) $, $K_\e$, and $u_\e(r_\e)$
are now as in (\ref{elements}), (\ref{kes1}), and (\ref{ta1}), respectively.
The following results hold true.

\bigskip
\noindent \textbf{Theorem 5.} \textit{Let $s\in\{1,\ldots,d\}$ be such that {\rm (\ref{ogran})} holds true.
Assume  that $\log d=o(\log\e^{-1})$ and that the quantity
$r_\e=r_\e(s)>0$ satisfies
$$\liminf_{\e\to 0}\frac{u_\e(r_\e)  }{\sqrt{2\log (d/s)}}>1.$$
Then as $\e\to 0$}
$$\sup_{s\in\mathcal{S}_d}\sup_{\eta\in{\mathcal H}_{d,s}}\sup_{\t\in\Theta_{\s,d}(r_\e)}s^{-1}\E_{\eta,\t}|\wh{\eta}(s_{\wh{m}})-\eta|\to 0.$$

\bigskip
\noindent
\textbf{Theorem 6.} \textit{Let $s\in\{1,\ldots,d\}$ be such that $s=o(d)$.
Assume that $\log d=o(\log\e^{-1})$ and that the quantity $r_\e=r_\e(s)>0$ satisfies
$$\limsup_{\e\to 0}\frac{u_\e(r_\e)  }{\sqrt{2\log (d/s)}}<1.$$
Then
$$\liminf_{\e\to 0}\inf_{\ti{\eta}}\sup_{\eta\in{\mathcal H}_{d,s}}\sup_{\t\in \Theta_{\s,d}(r_\e)} s^{-1}\E_{\eta,\t}|\eta-\ti{\eta}|> 0,$$
where the infimum is over all selectors $\ti{\eta}$ of a vector $\eta$ in model {\rm (\ref{model1})}.}

\bigskip
\noindent \textbf{Remark 2.}
We should remark that the best selection procedure yields exact variable selection only if the condition
$\liminf_{\e\to 0}\frac{u_\e(r_\e)  }{\sqrt{2\log d}+\sqrt{2\log s}}>1$ holds; at the same time, the best
selection procedure gives almost full variable selection if a milder condition $\liminf_{\e\to 0}\frac{u_\e(r_\e)  }{\sqrt{2\log (d/s)}}>1$
is met.

\section{Proofs of the Theorems}
In this section, we prove Theorems 3 and 4. The proofs of Theorems 5 and 6 go along the same lines and therefore are omitted.
Throughout the proof, the exponential bounds (\ref{expon1}) and (\ref{expon2}) on the tail probabilities of the statistics
$t_{j}(s)$ will be frequently used.

\bigskip\noindent \textbf{Proof of Theorem 3.} Let $m_0\in\{2,\ldots, M\}$ be such that
$$s_{m_0-1}\leq s< s_{m_0},$$
which implies that $s_{m_0}/s< d^{\Delta}$.
Then, using the definition of the selector $\hat{\eta}(s_{\hat{m}})$, we can write
\begin{gather}
\sup_{\eta\in{\mathcal H}_{d,s}}\sup_{\theta\in\Theta_{\s,d}(r_\e)}s^{-1}\E_{\eta,\t}|\hat{\eta}(s_{\hat{m}})-\eta|\nonumber\\ \leq
\sup_{\eta\in{\mathcal H}_{d,s}}\sup_{\theta\in\Theta_{\s,d}(r_\e)} s^{-1}\E_{\eta,\t}\left(|\hat{\eta}(s_{\hat{m}})-\eta| | \hat{m}< m_0\right)
\Pb_{\eta,\t}\left( \hat{m}< m_0 \right)\nonumber \\ +\sup_{\eta\in{\mathcal H}_{d,s}}\sup_{\theta\in\Theta_{\s,d}(r_\e)}s^{-1}\E_{\eta,\t}\left(|\hat{\eta}(s_{\hat{m}})-\eta| | \hat{m}\geq m_0\right)
\Pb_{\eta,\t}\left( \hat{m}\geq  m_0 \right)\nonumber\\
\leq \sup_{\eta\in{\mathcal H}_{d,s}}\sup_{\theta\in\Theta_{\s,d}(r_\e)} s^{-1}\E_{\eta,\t}\left(|\hat{\eta}(s_{\hat{m}})-\eta| | \hat{m}< m_0\right)
\Pb_{\eta,\t}\left( \hat{m}< m_0 \right)\nonumber\\
+ \sup_{\eta\in{\mathcal H}_{d,s}}\sup_{\theta\in\Theta_{\s,d}(r_\e)} (d/s)\Pb_{\eta,\t}\left( \hat{m}\geq m_0 \right)
=:I_1+I_2.
 \label{zvezda11}
\end{gather}
 To complete the proof, we need to show that $I_1$ and $I_2$ are both negligibly small when $\e$ is small.

Consider the term $I_1$ and observe that for all $\eta\in{\mathcal H}_{d,s}$ and $\theta\in\Theta_{\s,d}(r_\e)$,
\begin{gather*}
s^{-1}\E_{\eta,\t}\left(|\hat{\eta}(s_{\hat{m}})-\eta| | \hat{m}< m_0\right)\Pb_{\eta,\t}\left( \hat{m}< m_0 \right)\\  \leq s^{-1}\E_{\eta,\t}\left(|\hat{\eta}(s_{\hat{m}})-\hat{\eta}(s_{m_0})| \ \hat{m}< m_0\right)+
s^{-1}\E_{\eta,\t}\left(|\hat{\eta}(s_{m_0})-\eta| | \hat{m}< m_0\right)\Pb_{\eta,\t}\left( \hat{m}< m_0 \right)\\
\leq s^{-1}v_{m_0} + s^{-1}\E_{\eta,\t}|\hat{\eta}(s_{m_0})-\eta| ,
\end{gather*}
where by (\ref{deltacond}) and the choice of the sequences $\tau_d$ and $\Delta$
\begin{gather*}
s^{-1}{v_{m_0}}=\tau_d^{-1}(s_{m_0}/s)<{\tau_d}^{-1}d^{\Delta}=o(1).
\end{gather*}
Next, by definition of the set ${\cal H}_{d,s}$ of $s$-sparse $d$-dimensional  vectors $\eta$, we have
\begin{gather}
\sup_{\eta\in{\mathcal H}_{d,s}}\sup_{\theta\in\Theta_{\s,d}(r_\e)} s^{-1}\E_{\eta,\t}|\hat{\eta}(s_{m_0})-\eta| \leq
(d/s) \Pb_0\left(t_1(s_{m_0})>\sqrt{2\log(d/s_{m_0})+\delta \log d   }  \right)\nonumber\\ +
\sup_{\t_1\in\Theta_\s(r_\e)} \Pb_{\theta_1}\left(t_1(s_{m_0})\leq \sqrt{2\log(d/s_{m_0})+\delta \log d   }  \right)\label{nomer}
\end{gather}
where by (\ref{expon1}) the first summand in the above expression satisfies
\begin{gather*}
(d/s)\Pb_0\left(t_1(s_{m_0})>\sqrt{2\log (d/s_{m_0})+\delta\log d} \right) \\\leq (d/s)\exp\left(-\left(\log(d/s_{m_0})+(\delta/2)\log d\right)(1+o(1))  \right)\\
=O\left((s_{m_0}/s) d^{-\delta/2}\right)=O\left(d^{\Delta-\delta/2}\right)=o(1),
\end{gather*}
and the last equality is due to (\ref{deltacond0}) and (\ref{deltacond}).

To treat the second term on the right side of (\ref{nomer}), recall that $1<s_{m_0}/s<d^{\Delta}$. Then, by the assumption on the parameter $r_{\e}=r_{\e}(s)$ and the `continuity' of the function $u_\e(r_\e)$
as stated in (\ref{ucont}),
 using the fact that $\Delta \log d\to 0$ as $d \to \infty$,
 one can find a constant $\delta_1>0$ such that
for all sufficiently small $\e$
\begin{gather*}
r_\e\geq r^*_\e(s_{m_0})(1+\delta_1).
\end{gather*}
From this, using Proposition 4.1 in \cite{GI} and recalling formula (\ref{mo}),
\begin{eqnarray}\label{minus22}
\inf_{\t_1\in\Theta_{\s}(r_\e)} \Eb_{\t_1}\left( t_1(s_{m_0})\right)&\geq& \inf_{\t_1\in\Theta_{\s}(r^*_\e(s_{m_0})(1+\delta_1))} \Eb_{\t_1}\left( t_1(s_{m_0})\right)\nonumber\\
&\geq& (1+\delta_1)^2 \inf_{\t_1\in\Theta_{\s}(r^*_\e(s_{m_0}))} \Eb_{\t_1}\left( t_1(s_{m_0})\right)\geq (1+\delta_1)^2u_\e(r_\e^*(s_{m_0}))\nonumber\\
&=&(1+\delta_1)^2\sqrt{2\log (d/s_{m_0})}>\sqrt{2\log (d/s_{m_0})+\delta\log d },
\end{eqnarray}
where the last inequality follows from the fact that
$d^c\leq s_{m_0}<d^C$, which implies
$\delta\log d=o(\log (d/s_{m_0})).$
Thus as $\e\to 0$
\begin{gather}\label{beskon22}
\sqrt{2\log (d/s_{m_0})+\delta\log d }-\inf_{\t_1\in\Theta_{\s}(r_\e)} \Eb_{\t_1}\left( t_1(s_{{m_0}})\right)\to -\infty.
\end{gather}
Now (\ref{expon2}) in combination with (\ref{minus22}) and  (\ref{beskon22}) gives, uniformly in $\t_1\in\Theta_\s(r_\e)$,
\begin{gather*}
\Pb_{\t_1}\left(t_1(s_{m_0})\leq \sqrt{2\log (d/s_{m_0})+\delta\log d}   \right)\\
\leq \Pb_{\t_1}\left(t_1(s_{m_0})- \Eb_{\t_1}\left( t_1(s_{m_0})\right)\leq \sqrt{2\log (d/s_{m_0})+\delta\log d}-\inf_{\t_1\in\Theta_{\s}(r_\e)} \Eb_{\t_1}\left( t_1(s_{m_0})\right)   \right)\\
\leq \Pb_{\t_1}\left(t_1(s_{m_0})- \Eb_{\t_1}\left( t_1(s_{m_0})\right)\leq -\sqrt{2\log (d/s_{m_0})} \left[(1+\delta_1)^2-1+o(1)\right]\right)\\
\leq \exp\left( -\log (d/s_{m_0})\left[(1+\delta_1)^2-1+o(1)\right]^2 (1+o(1)) \right)\\
=O\left(  (s_{m_0}/d)^{\left[(1+\delta_1)^2-1\right]^2}\right)=o(1).
\end{gather*}
Putting everything together, we conclude that the first term on the right side of (\ref{zvezda11}) satisfies
\begin{gather}\label{i1}
I_1=o(1),\quad \e \to 0.
\end{gather}

Let us now show that
\begin{gather*}
I_2=\sup_{\eta\in{\mathcal H}_{d,s}}\sup_{\theta\in\Theta_{\s,d}(r_\e)} (d/s)\Pb_{\eta,\t}\left( \hat{m}\geq m_0 \right)=o(1).
\end{gather*}
By definition of $\hat{m}$,
for all $\eta\in{\mathcal H}_{d,s}$ and all $\theta\in\Theta_{\s,d}(r_\e)$,
\begin{gather*}
\Pb_{\eta,\t}\left( \hat{m}\geq m_0 \right)=\sum_{k=m_0}^{M}\Pb_{\eta,\t}\left( \hat{m}=k\right)\\=
\sum_{k=m_0}^M\Pb_{\eta,\t}\left(\exists\, i\in\{k,\ldots, M\}: |\hat{\eta}(s_{k-1})-\hat{\eta}(s_i)|>v_i\right)\\
\leq \sum_{k=m_0}^M\sum_{i=k}^M \Pb_{\eta,\t}\left(|\hat{\eta}(s_{k-1})-\hat{\eta}(s_i)|>v_i\right)\\
=\sum_{k=m_0}^M\sum_{i=k}^M \Pb_{\eta,\t}\left(\sum_{j=1}^d|\hat{\eta}_j(s_{k-1})-\hat{\eta}_j(s_i)|>v_i\right).
\end{gather*}
Now, we introduce independent events
\begin{gather*}
A_j(s)=\left\{t_j(s)>\sqrt{2\log (d/s)+\delta\log d}\right\},\quad j=1,\ldots,d,
\end{gather*}
and denote by $\overline{A_j(s)}$ the complement of $A_j(s)$. Observing that
for all $m_0\leq k\leq i\leq M$ the quantity $|\hat{\eta}_j(s_{k-1})-\hat{\eta}_j(s_i)|$
is non-zero only  if either $A_j(s_{k-1})\cap \overline{A_j(s_i)}$ or $ \overline{A_j(s_{k-1})}\cap {A_j(s_i)}$
occurs, we may continue
\begin{gather*}
\Pb_{\eta,\t}\left( \hat{m}\geq m_0 \right)\leq \sum_{k=m_0}^M \sum_{i=k}^M \Pb_{\eta,\t}\left( \sum_{j=1}^d\left[
\mathbb{I}\left(A_j(s_{k-1})\cap \overline{A_j(s_i)}  \right)+\mathbb{I}\left( \overline{A_j(s_{k-1})}\cap {A_j(s_i)} \right)\right]\right).
\end{gather*}
To bound this sum, we shall apply Bernstein's inequality  saying that if $\mathbb{X}_1,
\ldots, \mathbb{X}_d$ are independent random variables such that  for all  $j=1,\ldots,d$
and for some $H>0$
\begin{gather}\label{Bc}
\Eb(\mathbb{X}_j)=0\quad \mbox{and}\quad \left|\E (\mathbb{X}^m_j)\right|\leq \frac{\E(\mathbb{X}^2_j)}{2} H^{m-2} m!<
\infty,\quad m=2,3,\ldots,
\end{gather}
then  (see, for example, pp. 164--165 of \cite{Bernstein})
\begin{gather}\label{BI}
\max\left\{\Pb\left( \mathbb{S}_d\geq t\right), \Pb\left(\mathbb{S}_d \leq -t \right)\right\}\leq
\left\{\begin{array}{ll} \exp\left(-{t^2}/{4 B^2_d} \right)& \mbox{if}\quad0\leq t\leq B^2_d/H,\\
\exp\left(-{t}/{4H}  \right)& \mbox{if}\quad t\geq B^2_d/H,
\end{array}\right.
\end{gather}
where $\mathbb{S}_d=\sum_{j=1}^d \mathbb{X}_j$ and $B^2_d=\sum_{j=1}^d \Eb(\mathbb{X}^2_j).$ Observe that for independent
random variables $\mathbb{X}_1, \ldots, \mathbb{X}_d$ with the property
$$\Eb(\mathbb{X}_j)=0\quad \mbox{and}\quad\left|\mathbb{X}_j\right|\leq M,\quad  j=1,\ldots,d,$$
for some $M>0$,  the  Bernstein condition (\ref{Bc}) holds with  $H=M/3$. Below we will use Bernstein's inequality in the case of $t\geq B^2_d/H$.

To do this,  let us introduce random variables $\mathbb{X}_j=\mathbb{X}_j(s_{k-1},s_i)$,  $1\leq j\leq d,$ $ m_0 \leq k\leq M$, $k\leq i\leq M$, by the formula
\begin{eqnarray*}
\mathbb{X}_j&=&\mathbb{I}\left(A_j(s_{k-1})\cap \overline{A_j(s_i)}  \right)+\mathbb{I}\left( \overline{A_j(s_{k-1})}\cap {A_j(s_i)} \right)\\
 & &-\left[\Pb_{\eta,\t}\left(A_j(s_{k-1})\cap \overline{A_j(s_i)}  \right)+\Pb_{\eta,\t}\left( \overline{A_j(s_{k-1})}\cap {A_j(s_i)} \right) \right],
\end{eqnarray*}
and observe that $ |\mathbb{X}_j|\leq 4$, $j=1,\ldots,d$, and for all $\eta\in{\mathcal H}_{d,s}$ and $\theta\in\Theta_{\s,d}(r_\e)$
\begin{gather*}
\Eb_{\eta,\t}(\mathbb{X}_j)=0,\quad j=1,\ldots,d.
\end{gather*}
Before applying Bernstein's inequality, we show that for all $\eta\in{\mathcal H}_{d,s}$ and $\theta\in\Theta_{\s,d}(r_\e)$, and for all $m_0 \leq k\leq M$ and $k\leq i\leq M$
\begin{gather}\label{summa1}
\sum_{j=1}^d \left[\Pb_{\eta,\t}\left(A_j(s_{k-1})\cap \overline{A_j(s_i)}  \right)+\Pb_{\eta,\t}\left( \overline{A_j(s_{k-1})}\cap {A_j(s_i)} \right)\right] =o(v_i).
\end{gather}

We have
\begin{gather}
\sup_{\eta\in{\mathcal H}_{d,s}}\sup_{\theta\in\Theta_{\s,d}(r_\e)}\sum_{j=1}^d \left(\Pb_{\eta,\t}\left(A_j(s_{k-1})\cap \overline{A_j(s_i)}  \right)+\Pb_{\eta,\t}\left( \overline{A_j(s_{k-1})}\cap {A_j(s_i)} \right)\right)\nonumber\\
=(d-s) \left[\Pb_{0}\left(A_1(s_{k-1})\cap \overline{A_1(s_i)}  \right)+\Pb_{0}\left( \overline{A_1(s_{k-1})}\cap {A_1(s_i)} \right)\right] \nonumber\\
+ s \sup_{\theta_1\in\Theta_{\s}(r_\e)}\left[\Pb_{\t_1}\left(A_1(s_{k-1})\cap \overline{A_1(s_i)}  \right)+\Pb_{\t_1}\left( \overline{A_1(s_{k-1})}\cap {A_1(s_i)} \right)\right] \nonumber\\
\leq d\left[ \Pb_0\left(t_1(s_{k-1})>\sqrt{2\log (d/s_{k-1})+\delta\log d}  \right) +\Pb_0\left(t_1(s_{i})>\sqrt{2\log (d/s_{i})+\delta\log d} \right) \right]\nonumber\\
+ s\!\! \!\!\sup_{\theta_1\in\Theta_{\s}(r_\e)}\!\!\left[\Pb_{\t_1}\left(t_1(s_{k-1})\leq \sqrt{2\log (d/s_{k-1})+\delta\log d}   \right) \right.+ \left. \Pb_{\t_1}\left(t_1(s_{i})\leq \sqrt{2\log (d/s_{i})+\delta\log d}   \right)\right]\nonumber\\
=: J_1(s_{k-1},s_i)+ J_2(s_{k-1},s_i).\label{summa}
\end{gather}
Recalling (\ref{expon1}) and the relation $\tau_d d^{-\delta/2}\to 0$ as $d\to \infty$, we have
\begin{gather*}
d\Pb_0\left(t_1(s_{i})>\sqrt{2\log (d/s_{i})+\delta\log d} \right) \leq d\exp\left(-\left(\log(d/s_i)+(\delta/2)\log d\right)(1+o(1))  \right)
\\=O\left(s_i d^{-\delta/2}\right)= O\left(v_i\tau_d d^{-\delta/2}\right)=o(v_i).
\end{gather*}
Similarly, using the fact that $v_{k-1}<v_i$ when $k\leq i\leq M$, we obtain
\begin{gather*}
d\Pb_0\left(t_1(s_{k-1})>\sqrt{2\log (d/s_{k-1})+\delta\log d} \right) =o(v_{k-1})=o(v_i).
\end{gather*}
Therefore for all $m_0 \leq k\leq M$ and $k\leq i\leq M$
\begin{gather}\label{j1}
 J_1(s_{k-1}, s_i)=o(v_i).
\end{gather}

Consider the second term on the right side of  (\ref{summa}), $ J_2(s_{k-1},s_i)$.
First, note that for all $m_0 \leq k\leq M$ and $k\leq i\leq M$,
$$s<s_i\quad \mbox{and}\quad s<s_{k-1},\quad k\neq m_0.$$
and for $k=m_0$ one has $s_{k-1}=s_{m_0-1}\leq s$,
 which implies $s/s_{m_0-1}< d^{\Delta}$.
 Therefore, by the assumption on $r_{\e}=r_{\e}(s)$ and the `continuity' of the function $u_\e(r_\e)$ as cited in (\ref{ucont}),
 using the fact that $\Delta \log d\to 0$ as $d \to \infty$,
 one can find constants $\delta_2>0$ and $\delta_3>0$ such that
for all sufficiently small $\e$
\begin{gather*}
 r_\e\geq r^*_\e(s_{i})(1+\delta_2)\quad \mbox{and}\quad r_\e\geq r^*_\e(s_{k-1})(1+\delta_3)
\end{gather*}
when $ m_0 \leq k\leq M$ and $k\leq i\leq M$.
From this, for all sufficiently small $\e$, cf. (\ref{minus22}),
\begin{eqnarray}\label{minus}
\inf_{\t_1\in\Theta_{\s}(r_\e)} \Eb_{\t_1}\left( t_1(s_{i})\right)\geq
(1+\delta_2)^2\sqrt{2\log (d/s_i)}>\sqrt{2\log (d/s_i)+\delta\log d },
\end{eqnarray}
and hence as $\e\to 0$
\begin{gather}\label{beskon}
\sqrt{2\log (d/s_i)+\delta\log d }-\inf_{\t_1\in\Theta_{\s}(r_\e)} \Eb_{\t_1}\left( t_1(s_{i})\right)\to -\infty.
\end{gather}
It now follows from (\ref{expon2}), (\ref{minus}), and  (\ref{beskon}) that, uniformly in $\t_1\in\Theta_\s(r_\e)$,
\begin{gather*}
s\Pb_{\t_1}\left(t_1(s_{i})\leq \sqrt{2\log (d/s_{i})+\delta\log d}   \right)\\
\leq s\Pb_{\t_1}\left(t_1(s_{i})- \Eb_{\t_1}\left( t_1(s_{i})\right)\leq \sqrt{2\log (d/s_{i})+\delta\log d}-\inf_{\t_1\in\Theta_{\s}(r_\e)} \Eb_{\t_1}\left( t_1(s_{i})\right)   \right)\\
\leq s\Pb_{\t_1}\left(t_1(s_{i})- \Eb_{\t_1}\left( t_1(s_{i})\right)\leq -\sqrt{2\log (d/s_i)}\left[ (1+\delta_2)^2-1+o(1)\right]\right)\\
\leq s\exp\left( -\log (d/s_i)\left[(1+\delta_2)^2-1+o(1)\right]^2 (1+o(1)) \right)\\=
O\left( s (s_i/d)^{[(1+\delta_2)^2-1]^2}\right) =O\left( s_i (s_i/d)^{[(1+\delta_2)-1]^2}\right)=o(v_i).
\end{gather*}
Also, as relation (\ref{beskon}) continues to hold with $s_{k-1}$, $m_0\leq k\leq M$, instead of $s_i$, similar arguments yield
\begin{gather*}
s\Pb_{\t_1}\left(t_1(s_{k-1})\leq \sqrt{2\log (d/s_{k-1})+\delta\log d}   \right)=o(v_{k-1})=o(v_i),
\end{gather*}
which implies
\begin{gather}\label{j2}
J_2(s_{k-1},s_i)=o(v_i).
\end{gather}
Combining (\ref{summa}), (\ref{j1}) and (\ref{j2}), we arrive at (\ref{summa1}).
We see then by (\ref{summa1}) that
\begin{eqnarray*}
\sum_{j=1}^d \Eb_{\eta,\t}(\mathbb{X}_j^2)&=&\left(\sum_{j=1}^d\left[ \Pb_{\eta,\t}\left(A_j(s_{k-1})\cap \overline{A_j(s_i)}  \right)
+\Pb_{\eta,\t}\left( \overline{A_j(s_{k-1})}\cap {A_j(s_i)} \right)\right]\right)(1+o(1))\\&=&o(v_i).
\end{eqnarray*}
Therefore, the use of Bernstein's inequality as in (\ref{BI}) for the case of $t\geq B^2_d/H$ with $H=4/3$ gives
\begin{gather*}
I_2=\sup_{\eta\in{\mathcal H}_{d,s}}\sup_{\theta\in\Theta_{\s,d}(r_\e)} (d/s)\Pb_{\eta,\t}\left( \hat{m}\geq m_0 \right)\\
\leq \sup_{\eta\in{\mathcal H}_{d,s}}\sup_{\theta\in\Theta_{\s,d}(r_\e)} (d/s)\sum_{k=m_0}^M  \sum_{i=k}^M\Pb_{\eta,\t}
\left( \sum_{j=1}^d \mathbb{X}_j>v_i(1+o(1))\right)\\
\leq  (d/s)\sum_{k=m_0}^M  \sum_{i=k}^M
\exp\left(-({3v_i}/{16})(1+o(1))  \right) =O\left(M^2(d/s)\exp\left( -({3}/{16})v_{m_0} \right)\right)\\
=O\left(M^2(d/s)\exp\left(-({3d^{c}}/16\tau_d)   \right)  \right)=o(1).
\end{gather*}
This in combination with (\ref{zvezda11}) and (\ref{i1}) completes the proof of Theorem 3. \done

\subsection{Proof of Theorem 4} To prove the theorem, we first pick good prior distributions on $\eta=(\eta_j)$ and $\t=(\t_j)$.
Having done this, we bound the normalized minimax risk by the normalized Bayes risk
and show that the latter is strictly positive. The first part of the proof up to relation (\ref{lb}) go along the lines of
that of Theorem 2 in \cite{IngStep2014}, with $p=s/d$ instead of $p=d^{-\beta}$.

Let
$\t^*_j=(\t^*_{j,k})_{k\in\mathbb{Z}}$ be the extremal sequence in the problem (the same for all $j=1,\ldots,d$)
\begin{gather*}
\frac{1}{2\e^4}\sum_{k\in\mathbb{Z}} \t_{j,k}^4 \to \inf_{\t_j\in \Theta_{\s}(r_\e)}.
\end{gather*}
Let the prior distribution of a `vector'  $\t=(\t_1,\ldots,\t_d)\in\Theta_{\s,d}(r_\e)$ be of the form
$$\pi_{\theta}(d\t)=\prod_{j=1}^d \pi_{\t_j}(d\t_j), \quad
\pi_{\t_j}(d\t_j)=\prod_{1\leq|k|\leq K_\e} \left(\frac{\delta_{-\t^*_{j,k}}+\delta_{\t^*_{j,k}}}{2}\right)(d\t_{j,k}),$$
where $\delta_x$ is the $\delta$-measure that puts a pointmass 1 at $x$.
Denote by
$$ p=s/d$$
the portion of non-zero components of vector $\eta=(\eta_1,\ldots,\eta_d)\in{\mathcal H}_{d,s}$.
The prior distribution of $\eta$ is naturally defined to be
$$\pi_\eta(d\eta)=\prod_{j=1}^d \pi_{\eta_j}(d\eta_j), \quad \pi_{\eta_j}(d\eta_j)=\left((1-p)\delta_0+p\delta_1\right)(d\eta_j).$$
Then, assuming that $\t=(\t_j)$ and $\eta=(\eta_j)$ are independent, we get
\begin{eqnarray*}
R_\e&:=&\inf_{\ti{\eta}}\sup_{\eta\in{\mathcal H}_{d,s}}\sup_{\t\in \Theta_{\s,d}(r_\e)} s^{-1}\E_{\eta,\t}|\eta-\ti{\eta}|
\geq  s^{-1}\inf_{\ti{\eta}}\E_{\pi_\eta} \E_{\pi_\t} \E_{\eta,\t}|\eta-\ti{\eta}|\\
&=& s^{-1}\inf_{\ti{\eta}} \E_{\pi_\eta} \E_{\pi_\t} \E_{\eta,\t}\sum_{j=1}^d |\eta_j-\ti{\eta}_j|
= s^{-1}\inf_{\ti{\eta}}\sum_{j=1}^d \E_{\pi_{\eta_j}} \E_{\pi_{\t_j}} \E_{\eta_j\t_j}|\eta_j-\ti{\eta}_j|,
\end{eqnarray*}
where the infimum is over all selectors $\ti{\eta}=(\ti{\eta}_j)$ and $\E_{\eta_j\t_j}$
is the expected value that corresponds to the measure $\Pb_{\eta_j \t_j}$ induced by the observation $X_j=(X_{j,k})_{1\leq|k|\leq K_\e}$
consisting of independent random variables $X_{j,k}$ that follow normal distributions $N(\eta_j \t_{j,k},\e^2)$.

Consider the mixture of distributions
given by the formula
\begin{multline}\label{density}
\Pb_{\pi,\eta_j}(d X_j)=\E_{\pi_{\t_j}} \Pb_{\eta_j \t_j} (dX_{j,k})
= \prod_{1\leq|k|\leq K_\e}  \left(\frac{N(-\eta_j \t^*_{j,k},\e^2)+N(\eta_j \t^*_{j,k},\e^2)}{2}\right)(d X_{j,k}).
\end{multline}
In particular, when $\eta_j=0$,
$\Pb_{\pi,0}(d X_j)=\prod\limits_{1\leq|k|\leq K_\e} N(0,\e^2) (d X_{j,k})$.
Using the notation
\begin{gather}\label{v}
v^*_{j,k}= \frac{\t^*_{j,k}}{\e},
\end{gather}
we obtain with respect to the probability measure $\Pb_{\pi,\eta_j}$
$$
Y_{j,k}:=\frac{X_{j,k}}{\e}=\eta_jv^*_{j,k}+\xi_{j,k} \stackrel{\rm ind.}\sim{{N}}(\eta_j v^*_{j,k},1),\quad 1\leq j\leq d, \quad 1\leq|k|\leq K_\e.
$$
Next, denoting $Y_j=(Y_{j,k})_{1\leq|k|\leq K_\e}$, we may rewrite the likelihood ratio in the form
\begin{gather}\label{lr}
\frac{d\Pb_{\pi,\eta_j}}{d\Pb_{\pi,0}}(Y_j)
= \prod_{1\leq|k|\leq K_\e}\exp\left( -\frac{\eta_j (v^*_{j,k})^2}{2} \right) \cosh \left( {\eta_j v^*_{j,k}} Y_{j,k} \right).
\end{gather}
From this, using the fact that each $\eta_j$ takes on only two values, zero and one, with respective probabilities
$(1-p)$ and $p$, we may continue
\begin{multline}\label{rep}
R_\e\geq s^{-1}\sum_{j=1}^d \inf_{\ti{\eta}_j}\E_{\pi_{\eta_j}}  \E_{\pi,\eta_j} |\eta_j-\ti{\eta}_j|
=
s^{-1}\sum_{j=1}^d \inf_{\ti{\eta}_j}\left[(1-p)\E_{\pi,0}(\ti{\eta}_j)+p \E_{\pi,1}(1-\ti{\eta}_j)\right],
\end{multline}
where $\inf_{\ti{\eta}_j} (1-p)\E_{\pi,0}(\ti{\eta}_j)+p \E_{\pi,1}(1-\ti{\eta}_j)$
is the Bayes risk in the problem of testing two simple hypotheses
$$
H_0: \Pb=\Pb_{\pi,0}\quad\mbox{vs.}\quad H_1: \Pb=\Pb_{\pi,1},
$$
with the probability measures $\Pb_{\pi,0}$ and $\Pb_{\pi,1}$
defined according to (\ref{density}).
In particular, under the null hypothesis, the vector
 $Y_j=(Y_{j,k})_{1\leq|k|\leq K_\e}$ has a normal distribution with
density function
$p_{\pi,0}(t)\!=\!\prod\limits_{1\leq|k|\leq K_\e}(2\pi)^{-1/2}\exp(\! -t_k^2/2)$,
$t=(t_k)_{1\leq|k|\leq K_\e}$.
By (\ref{lr}) the likelihood ratio in this problem becomes
$$\Lambda_{\pi}=\Lambda_{\pi}(Y_j)=\frac{d\Pb_{\pi,1}}{d\Pb_{\pi,0}}(Y_j)=
\prod_{1\leq|k|\leq K_\e}\exp\left( -\frac{ (v^*_{j,k})^2}{2} \right) \cosh \left( { v^*_{j,k}} Y_{j,k} \right),   $$
and the optimal (Bayes) test $\eta_B$ that minimizes the Bayes risk in hand has the form (see, for example, \cite[Sec.
\,8.11]{DeGr})
$$\eta_B(Y_j)=\mathbb{I}\left(\Lambda_{\pi}(Y_j)\geq \frac{1-p}{p}\right).$$
Using this, we infer from (\ref{rep}) that
\begin{eqnarray}\label{lb}
R_\e&=&\inf_{\ti{\eta}}\sup_{\eta\in{\mathcal H}_{d,s}}\sup_{\theta\in\Theta_{\s,d}(r_\e)}s^{-1}\Eb_{\eta,\t}|\eta-\ti{\eta}|\nonumber \\
&\geq& (d/s)\Pb_{\pi,0}\left(\Lambda_{\pi}(Y_1)\geq \frac{1-p}{p}\right)+\Pb_{\pi,1}\left(\Lambda_{\pi}(Y_1)< \frac{1-p}{p}\right)=:A_\e+B_\e.
\end{eqnarray}
where, under the $\Pb_{\pi,\eta_1}$-probability with $\eta_1\in\{0,1\}$, the vector $Y_1=(Y_{1,k})_{1\leq |k|\leq K_\e}$ has independent normal components
$$Y_{1,k}=\eta_1v^*_{1,k}+\xi_{1,k}\sim{{N}}(\eta_1 v^*_{1,k},1),\quad  1\leq |k|\leq K_\e.$$
It now follows from (\ref{lb}) that the minimax risk $R_\e$ is positive if at least one of the terms, $A_\e$ or $B_\e$, is positive.
Let us prove that for all sufficiently small $\e$ the probability $B_\e$ is separated from zero.

Recall that $d=d_\e\to \infty$ and $s=s_d=o(d)$ as $\e\to 0$. Put $$H=H_\e=\log\left(\frac{1-p}{p}\right)\sim \log (d/s),$$
and introduce the random variable $$\lambda_{\pi}=\lambda_{\pi}(Y_1):=\log \Lambda_{\pi}(Y_1).$$
Using the notation $\Pb_0$ for $\Pb_{\pi,0}$, consider the
probability measure $\Pb_h$, depending on a positive parameter $h=h_\e,$ that is defined by the formula
\begin{gather*}
\frac{d\Pb_h}{d \Pb_0}(Y_1) :=\frac{\exp({h\lambda_{\pi}(Y_1)})}{\Psi(h)},\quad \Psi(h)=\Eb_{\Pb_0}\exp({h\lambda_{\pi}(Y_1)}).
\end{gather*}
With the parameter $h>0$ chosen to satisfy
$$\Eb_{\Pb_h}\lambda_{\pi}=H,$$  we have
(see Lemma 2  in \cite{IngStep2014})
\begin{gather}\label{hs}
h\sim \frac12+\frac{H}{u_\e^2}=O(1),
\end{gather}
and (see formula (45) in \cite{IngStep2014})
\begin{gather}\label{psi}
\Psi(h)=\exp\left(\frac{h^2-h}{2}\,u_\e^2(1+o(1))   \right),
\end{gather}
where for notational simplicity we use $u^2_\e$ for $u^2_\e(r_\e)$.

 We have
\begin{eqnarray}
B_\e&=&\Eb_{\pi,1}\left(\mathbb{I}\left(\lambda_{\pi}(Y_1)< H\right)\right)=
\Eb_{\pi,0}\left(\exp(\lambda_{\pi}(Y_1))
\mathbb{I}\left(\lambda_{\pi}(Y_1))<H \right)  \right)\nonumber \\&=&
\Eb_{h}\left(\frac{d\Pb_0}{d \Pb_h}(Y_1) \exp(\lambda_{\pi}(Y_1))
\mathbb{I}\left(\lambda_{\pi}(Y_1))<H \right) \right)\nonumber\\&=&
\Psi(h)\Eb_{h}\left(\exp[(1-h)\lambda_{\pi}(Y_1)]
\mathbb{I}\left(\lambda_{\pi}(Y_1))<H \right) \right).\label{Be}
\end{eqnarray}
By Lemma 3 in \cite{IngStep2014}, the standardized random variable
$$Z_h:=\frac{\lambda_{\pi}-\mu_h}{\s_h}, $$
where
$$\quad \mu_h=\Eb_{\Pb_h}(\lambda_{\pi})= u_\e^2(h-1/2)(1+o(1)),\quad \s^2_h=\Var_{\Pb_h}(\lambda_{\pi})= u_\e^2(1+o(1)),$$
converges in $\Pb_h$-distribution to an $N(0,1)$. Therefore the statistic $\lambda_{\pi}(Y_1)$ on the right side of (\ref{Be})
is nearly a normal $N(H,u_\e^2)$ random variable.

Next, by assumption and the `continuity' of $u_\e$ as stated in (\ref{ucont}), for some constant $\delta_1>0$
\begin{gather*}%\label{usl1}
u_\e/\sqrt{\log(d/s)}\leq \sqrt{2}(1-\delta_1),
\end{gather*}
provided $\e$ is small enough.
%for all sufficiently small $\e$
This and formula (\ref{hs}) give the inequality $1-h<0$, which implies for all $y\in\mathbb{R}^{2K_\e}$ and all sufficiently small $\e$
\begin{gather*}
\exp[(1-h)\lambda_{\pi}(y)]
\mathbb{I}\left(\lambda_{\pi}(y))<H \right)\leq \exp\left[(1-h)H \right]\sim (d/s)^{1-h}\leq {\rm const}.
\end{gather*}
Then, by the dominant convergence theorem, the replacement of $\lambda_{\pi}(Y_1)$ by an $N(H,u_\e^2)$ on the right side (\ref{Be}) and the use
of  (\ref{hs}) and (\ref{psi}) yield for all sufficiently small $\e$
\begin{gather*}
B_\e\sim \exp\left(\frac{h^2-h}{2}u_\e^2 \right)\int_{-\infty}^{H}\exp\left[(1-h)x \right]
\frac{1}{\sqrt{2\pi}u_\e} \exp\left( -\frac{(x-H)^2}{2u_\e^2} \right) \, dx\\
=\exp\left(\frac{h^2-h}{2}u_\e^2 +H(1-h)+\frac{(1-h)^2 u_\e^2}{2}\right)\int_{-\infty}^H \frac{1}{\sqrt{2\pi}u_\e}
\exp\left( -\frac{\left(x-(H+(1-h)u_\e^2)\right)^2}{2u_\e^2} \right)\, dx\\
\sim\exp(0)\int_{-\infty}^H \frac{1}{\sqrt{2\pi}u_\e}
\exp\left( -\frac{\left(x-(H+(1-h)u_\e^2)\right)^2}{2u_\e^2} \right)\, dx\\ \geq
\int_{-\infty}^{H+(1-h)u_\e^2} \frac{1}{\sqrt{2\pi}u_\e}
\exp\left( -\frac{\left(x-(H+(1-h)u_\e^2)\right)^2}{2u_\e^2} \right)\, dx=1/2.
\end{gather*}
From this
$$\liminf_{\e\to 0}R_\e\geq \liminf_{\e\to 0}B_\e\geq 1/2>0,$$
and the proof of Theorem 4 is complete. \done

\section{Concluding remarks}
In the context of variable selection in high dimensions, in both regression and white noise settings, simple thresholding
provides plausible alternative to the lasso for a large range of problems.
As a statistical tool, thresholding strategy is simple in nature and is not as computationally demanding as the lasso,
especially in very high dimensional problems.
At the same time, it is capable of doing at least as good as the lasso, or even better (see our Theorems 1 to 6,
Theorems 9 to 11 in \cite{GJW2009}, and Theorems 1 and 2 in \cite{IngStep2014} for details).
In light of these facts, we support the viewpoint of Genovese et al. \cite{GJW2009}
that for sparse high-dimensional regression problems a simple thresholding procedure merits further investigation.

To conclude our study, we point out possible directions for extending the results obtained in this paper. % and earlier in work \cite{IngStep2014}.
For the two function spaces ${\cal F}_\s$ at hand, it might be of interest to produce asymptotically exact and almost full selectors
 in very high dimensional settings when the conditions $\log d=o\left(\e^{-2/(2\s+1)}  \right)$
and $\log d=o(\log \e^{-1})$ on the growth of $d$ as a function of $\e$ are violated.

The setup of inverse problems, where the observations are $X_\varepsilon = Kf + \varepsilon W$, with $K$ being
a linear operator such that $K^\star K$ is compact, translates into a Gaussian sequence model with 
heterogenous observations $X_{j,k} = \eta_j \theta_{j,k} + \epsilon v_k \xi_{j,k}$, where $v_k^{-2}$ are
the eigenvalues of $K^\star K$. This case, which extends our setup, can be treated by using  
the sharp testing results for the inverse problems obtained in \cite{ISS}.

Furthermore, handling the problem of variable selection in a sequence space model, general ellipsoids
$\{\t\in l_2(\mathbb{Z}): \sum_{k\in\mathbb{Z}} c_k^2\t_k^2\leq 1\}$ in $l_2(\mathbb{Z})$,
with semi-axes $c_k$ decreasing fast enough, could be studied.
A more complicated model, in which a $d$-variate regression function $f$ admits a decomposition to a sum of $k$-variate components, with $k\geq 2$
and only a small number $s$ of these components being non-zero,
also deserves some attention. 

Eliminating the assumption of \textit{known} parameter $\s$ leads to the problem of adapting the
proposed selection procedures to the possible values of $\s$.

To pursue more practical goals, one can try to translate the results obtained for an
additive $s$-sparse Gaussian white noise model to the corresponding discrete regression model for which 
the corresponding detection problem was solved in \cite{Abram}.

\end{document}